\documentclass[graybox]{svmult}

\usepackage{mathptmx}       
\usepackage{helvet}         
\usepackage{courier}        
\usepackage{type1cm}        
\usepackage{makeidx}         
\usepackage{graphicx}        
\usepackage{multicol}        
\usepackage[bottom]{footmisc}

\makeindex             

\begin{document}

\title*{Stochastic geodesics}
\author{Ana Bela Cruzeiro and Jean-Claude Zambrini}
\institute{Ana Bela Cruzeiro \at Grupo de F\'{i}sica-Matem\'atica and Dep. Matem\'atica I.S.T., Universidade de Lisboa, Av. Rovisco Pais, 1049-001 Lisboa, Portugal \email{ana.cruzeiro@tecnico.ulisboa.pt} \and
Jean-Claude Zambrini \at Grupo de F\'{i}sica-Matem\'atica, Fac. Ci\^encias Universidade de Lisboa, Campo Grande Ed. C6, 1749-016 Lisboa, Portugal \email{jczambrini@fc.ul.pt}}

\maketitle

\texttt{We describe, in an intrinsic way and using the global chart provided by It\^o's parallel transport,  a generalisation of the notion of geodesic (as critical path of an energy functional) to diffusion processes on Riemannian manifolds. These stochastic processes are no longer smooth paths but
they are still critical points of a regularised stochastic energy functional. We consider stochastic geodesics on compact Riemannian manifolds and also on (possibly infinite dimensional)  Lie groups.
Finally the question of existence of such  stochastic geodesics is  discussed: we show how it can be  approached via forward-backward stochastic differential equations.}

\section{Introduction}

The notion of geodesic in Riemannian manifolds appeared first in a lecture of Riemann, in 1854. Originally, it was referring to the shortest path between two points on Earth's surface. Nowadays, given an affine connection like the one of Levi-Civita, it can also be defined as a curve whose tangent vectors remain parallel when transported along the curve. In Theoretical Physics it is in General Relativity that this notion played a key r\^ole.

In a stochastic framework, a generalisation of geodesic curve is described. It corresponds to a critical path for some generalised action functional. The concept is reminiscent of Feynman path
integral approach to Quantum Mechanics (\cite{AKKM2008}) but for well defined probability measures on path spaces. It involves, in particular, regularisation of the second order in time classical dynamical equations, which is not traditional in Stochastic Analysis.

The derived equations of motion are of Burgers type. When considering flows which keep the volume measure invariant one obtains Navier-Stokes equations. This point of view was developed in
in \cite{CiC2007}, \cite{AC2012}  and \cite{ACC2014} in particular.  It is currently being investigated (c.f  \cite{CCR2018} as well as \cite{C2020} for a review on this subject).

After a short survey of classical geodesics on Riemannian manifolds, Cartan's frame bundle approach and its relation with the horizontal and Laplace-Beltrami operators are recalled.

Stochastic Analysis of diffusions on manifolds along the line of It\^o-Ikeda-Watanabe is given, together with It\^o's associated notion of parallel transport. Then one comes back to one of the historic definitions of geodesics, namely as critical points of an Action functional. The regularisations associated with the critical diffusion provide the appropriate generalised energy functional. The same strategy applies to geodesics on Lie groups.

It is also shown how, if needed, stochastic geodesics can be characterised via stochastic forward-backward SDEs.

It is a special pleasure to dedicate this paper to Sergio Albeverio as a modest sign of recognition for his faithful friendship along the years.

\section{Geodesics on Riemannian manifolds}
\label{sec:2}

We shall denote by $M$  a $d$-dimensional compact Riemannian manifold and $g$ its metric tensor.  Given $m\in M$, if $u, v$ are vectors in the tangent space $T_m (M)$ the Riemannian inner product is given in local chart by
\begin{eqnarray*}
g_m (u,v)=(g_{i,j} u^i v^j )(m)
\end{eqnarray*}
Here and in the rest of the paper we adopt Einstein summation convention.

The Levi-Civita covariant derivative of a vector field $z$ has the expression
\begin{eqnarray*}
[\nabla_k z]^j=\frac{\partial}{\partial m^k} z^j +\Gamma_{k,l}^j z^l,
\end{eqnarray*}
where $\Gamma$ denotes the corresponding Christoffel symbols in the local chart; explicitly,
\begin{eqnarray}\label{eq1}
\Gamma_{k,l}^j =\frac{1}{2} \left( \frac{\partial}{\partial m^k} g_{i,l}
+\frac{\partial}{\partial m^l} g_{k,i}-\frac{\partial}{\partial m^i} g_{k,l}\right) g^{j,i}
 \end{eqnarray}

Given a smooth curve $t\rightarrow \varphi (t)\in M$, the parallel transport of a vector field
$z$  along this curve is defined by the condition of zero covariant derivative of $z$ in the $\dot \varphi$
direction,

\begin{eqnarray}\label{eq2}
\nabla_{\dot \varphi (t)} z(t) =0~~\hbox{or}~~\dot z^j =-\Gamma_{k,l}^j \dot \varphi^k z^l.
\end{eqnarray}
 Its solution,
$z(t)=t^{\varphi}_{t\leftarrow 0}(z(0))$, the parallel transport of $z$ along the curve,
 is an Euclidean isomorphism between tangent spaces:

\begin{eqnarray*}
t^{\varphi}_{t\leftarrow 0}:T_{\varphi (0)} (M)\rightarrow T_{\varphi (t)}(M).
\end{eqnarray*}

Consider  the curves minimising the length
\begin{eqnarray*}
{\cal J}(\gamma )= \int l(\gamma ,\dot \gamma) dt, ~~~~~l(\gamma ,\dot \gamma ) =
\sqrt{g_{i,j}\dot \gamma^i \dot
\gamma^j}
\end{eqnarray*}
and therefore satisfy the Euler variational equation
\begin{eqnarray*}
\frac{d}{dt} (\frac{g_{i,j} \dot \gamma^j}{l})=\frac{1}{2l} \partial_i
(g_{j,k})\dot \gamma^j \dot \gamma^k .
\end{eqnarray*}
Replacing $dt$ by $ds$ (where $s$ is the arc lenght) we obtain
\begin{eqnarray*}
\frac{d}{dt} (g_{i,j} \dot \gamma^j \frac{1}{ds})-\frac{1}{2}\frac{1}{ds}
\partial_i (g_{j,k})\dot \gamma^j \dot \gamma^k =0
\end{eqnarray*}
and also
\begin{eqnarray*}
g_{i,j} \frac{d^2 \gamma^j}{ds^2}+\partial_k (g_{i,j})\frac{d\gamma^k }{ds}
\frac{d\gamma^j}{ds} -\frac{1}{2} \partial_i (g_{j,k})\frac{d\gamma^k }{ds}
\frac{d\gamma^j}{ds} =0
\end{eqnarray*}
Multiplying both members by $g^{\alpha ,i}$ we obtain the following classical
form of the geodesic equation:
\begin{eqnarray}\label{eq3}
\frac{d^2 \gamma^\alpha}{ds^2} +\Gamma_{j,k}^\alpha \frac{d\gamma^j}{ds}
\frac{d\gamma^k}{ds}=0 
\end{eqnarray}
or $\nabla_{\dot \gamma} \dot \gamma =0$.

A curve satisfying the last equation is called a \it geodesic \rm for the corresponding Riemannian metric. It is also well known that geodesics (defined in a time interval $[0,T]$)  are characterised as being critical paths of the (kinetic) energy functional
\begin{eqnarray}\label{eq4}
{\cal E} (\gamma )=\int_0^T ||\dot \gamma (t)||^2 dt=\int_0^T g_{i,j}(\gamma (t) )\dot \gamma^i (t)\dot \gamma^j (t) dt
\end{eqnarray}
By critical it is meant that, for every family of smooth curves (variations of $\gamma$)
$\gamma_{\epsilon}$ starting (at time $0$) and ending (at time $T$) at $\gamma (0)$ and $\gamma (T)$ resp.,  we have
$\frac{d}{d\epsilon}|_{\epsilon =0} {\cal E}(\gamma_{\epsilon})=0$.

\section{The frame bundle and the Laplacians}
\label{sec:3}
The bundle of orthonormal frames over $M$ is defined by
\begin{eqnarray*}
O(M)=\{ (m,r): m\in M, r: R^d \rightarrow T_m (M) ~~\hbox{is an Euclidean isometry}\}
\end{eqnarray*}
The map $\pi :O(M)\rightarrow M,~\pi (m,r)=m$ is the canonical projection.

Let $e_i ,~i=1,...d$ denote the vectors of the canonical basis of $R^d$ and
$\gamma_i$ denote the (unique) geodesic on $M$ such that $\gamma_i (0)=m,~
\left.\frac{d}{dt}\right|_{t=0} \gamma_i (t) =r(e_i )$. Let $(\gamma_i (t), r_i (t))$ represent the parallel transport of $r$ along $\gamma_i,~
\nabla_{\dot \gamma_i} r_i =0,~~~r_i (0)=\hbox{Id}$.
Then
\begin{eqnarray*}
A_i (m,r) =\left.\frac{d}{dt}\right|_{t=0} r_i (t)
\end{eqnarray*}
are called the horizontal vector fields on $M$.

Denote by $\Theta$ the one-form defined on $O(M)$ with values
in $ R^d \times \hbox{so}(d)$ such that  $<\Theta ,A_i >=(e_i ,0)$; $\Theta =(\theta, \omega )$,
with $\omega (m,r)=r^{-1} dr$ the Maurer-Cartan form on the orthogonal group $O(d)$. Its structure equations are given by
\begin{eqnarray*}
\left\{\begin{array}{lll}
d\theta =\omega \wedge \theta \\
d\omega = \omega \wedge \omega +\Omega (\theta \wedge \theta ),
\end{array}\right.
\end{eqnarray*}
where $\Omega$ denotes the curvature tensor:
\begin{eqnarray*}
\Omega (X,Y,Z)=(\nabla_X \nabla_Y -\nabla_Y \nabla_X -\nabla_{[X,Y]} )Z,
\end{eqnarray*}
where $[X,Y]$ denotes the bracket of two vector fields.
Recall also that the Ricci tensor ($\hbox{Ricci}_{kl}$) is the trace of the curvature, taken in the second and third entries.

In particular $\theta (A_k )=e_k$ and $\omega (A_k )=0$.
The horizontal Laplacian on $O(M)$ is the second order differential operator
\begin{eqnarray}\label{eq5}
\Delta_{O(M)} =\sum_{k=1}^d {\mathcal L}_{A_k}^2
\end{eqnarray}
where $ {\mathcal L}_{A_k}$ denotes the Lie derivative along the vector field $A_k$.
For every smooth function $f$ defined on $M$ we have
\begin{eqnarray*}
\Delta_{O(M)} (f \circ \pi )=(\Delta_M f) \circ \pi
\end{eqnarray*}
where $\Delta_M $ is the Laplace-Beltrami operator on $M$. This operator is
expressed in local  coordinates by
\begin{eqnarray}\label{eq6}
\Delta_M f= g^{i,j}[\frac{\partial^2 f}{\partial m^i \partial m^j}-\Gamma_{i,j}^k 
\frac{\partial f}{\partial m^k}].
\end{eqnarray}

\section{Stochastic Analysis on  manifolds}
\label{sec:4}
We are going to consider stochastic diffusions associated to elliptic operators on $M$ of the form
\begin{eqnarray}\label{eq7}
L_u f := \frac{1}{2} \Delta_M f +\partial_u f
\end{eqnarray}
in the sense of It\^o stochastic calculus. Here $u$ is a  possibly time-dependent,  smooth (at least $C^2$) vector field on $M$. In local coordinates the diffusion with generator $L_u$ can be written as
\begin{eqnarray}\label{eq8}
dm^j (t)=\sigma_k^j  dx^k (t) -\left( \frac{1}{2}g^{m,n} \Gamma^j_{m,n} -u^j \right) dt
\end{eqnarray}
where $\sigma =\sqrt g$ and $x_k $ are independent real-valued Brownian motions.

We consider the horizontal lift of these $M$-valued diffusion processes.
Denote by $u_k$ the functions defined on $O(M)$ by
\begin{eqnarray*}
U_k (r)=r(e_k ). u_{\pi (r )}.
\end{eqnarray*}
Then $\tilde u =\sum_k U_k A_k $ satisfies $\pi^{\prime} (\tilde u)=u$ ($\pi^{\prime}$ being the derivative of the canonical projection $\pi$).

Denoting by $x$ a sample path of the standard Brownian motion on $R^d$, $x(t), t\in [0,T],~ x(0)=0$, we consider the following Stratonovich stochastic differential equation  on $O(M)$:
\begin{eqnarray}\label{eq9}
dr_x (t)=\sum_{k=1}^d A_k (\circ dx^k (t)+U_k dt),~~~
r_x (0)=0
\end{eqnarray}
with $\pi (r_0 )=m_0$.
In local coordinates $(m^i , e_{\alpha}^i )$ on $O(M)$ and if $r(t)=(m(t),e(t))$ we have,

\begin{eqnarray*}
\left\{\begin{array}{lll}
dm^i (t) =e_{\alpha}^i \circ (dx^{\alpha} (t) +u^{\alpha}dt)\\
de_{\alpha}^i (t)= -\Gamma_{j,k}^i (m(t)) e_{\alpha}^k (t)\circ dm^j (t),
\end{array}\right.
\end{eqnarray*}

If $a\in M$,  we  denote the  path space of the manifold-valued paths starting from $a$ by
\begin{eqnarray*}
P_{a} (M)=\{p:[0,T]\rightarrow M, p(0)=a, p~~ \hbox{continuous}\}.
\end{eqnarray*}
The diffusion $m(t)$ has for generator the operator $L_u$. We refer to \cite{IkedaWatanabe1989} for a detailed exposition of diffusions
on Riemannian manifolds constructed on the frame bundle.

For each vector field $u$ the operator $L_u $ and the operator ${\mathcal L}_U =\frac{1}{2} \Delta_{O(M)}+\partial_U$ induce on the path spaces $P_{m_0} (M)$ and $P_{m_0} (O(M))$, respectively, two probability measures, namely the laws of the corresponding diffusion processes. The projection map $\pi$ realizes an isomorphism between these two probability spaces. 

Let  the path space $P_{0} (R^d )$ be endowed with the law of the process
$dy (t)= (\circ dx (t) +U) (y(t))$, $t\in [0,T]$ and $P_{m_0} (M)$ with the law of the diffusion $p$ with generator 
$L_u$).
Consider the It\^o map $I:P_{0} (R^d )\rightarrow P_{m_0} (M)$ defined by
\begin{eqnarray*}
I(x)(t)=\pi (r_x (t))
\end{eqnarray*}

This map is a.s. bijective and provides an isomorphism between the  corresponding probability measures (\cite{Malliavin1995}).

Even though $p$ is not differentiable in time, It\^o has shown that one can still define a parallel transport along $p$, which is the isomorphism from $T_{p(s)} (M)\rightarrow T_{p(t)} (M)$ given by
\begin{eqnarray*}
t^p_{t\leftarrow s}:= r_x (t) r_x (s)^{-1}.
\end{eqnarray*}

The differentiability of $r_x (t)$ with respect to variations of the Brownian motion $x$ was
studied in \cite{FangMalliavin1993} and \cite{ABCMalliavin1996} within the framework of Malliavin Calculus \cite{M1997}, \cite{B1984} (c.f. also \cite{ABCMalliavin1998} for the case of the Brownian motion with drift). 

Denote ${\mathcal D}_\alpha^\beta ={\mathcal L}_{A_\alpha} u^\beta$. The following result holds:

\vskip 5mm

\begin{proposition}\label{prop1} Given a process of bounded variation in time  $h:P_{0}(R^d)\times [0,T]\rightarrow R^d$, we have, using the notations of section $3$,
\begin{eqnarray}\label{eq10}
<\theta, \frac{d}{d\epsilon}_{|_{\epsilon =0}} r_{x+\epsilon h}>=\zeta,~~~~
<\omega, \frac{d}{d\epsilon}_{|_{\epsilon =0}} r_{x+\epsilon h}>=\rho
\end{eqnarray}
where $\zeta$ and $\rho$ are determined by  It\^o (and Stratonovich) stochastic differential equations
\begin{eqnarray}\label{eq11}
d\zeta (t)=\dot h (t) dt\ -[\frac{1}{2}\hbox{Ricci} +{\mathcal D}] (h (t)) dt -\rho (t)dx(t)\\
d\rho =\Omega (\circ dx +udt, h )\nonumber
\end{eqnarray}
with initial conditons $\zeta (0)=0$, $\rho (0)=0$.
\end{proposition}

The result above is still valid for pinned Brownian motion, namely when $p(T)$ is fixed. Then the variations are equal to zero not only at the initial but
also at this final time.
The  sigma-algebra and filtration on the corresponding path space are the usual ones,
generated by the coordinate maps and generated by the coordinate maps up to time t,
respectively. We refer to \cite{D1994} for more details.

\section{Stochastic geodesics}
\label{sec:5}

We shall consider stochastic geodesics as processes which are critical points of some energy functional generalising the classical deterministic one. 
Since the stochastic processes,
diffusions on the manifold, are no longer differentiable in time, some notion of generalised velocity has
to replace the usual time derivative.

If $\xi (\cdot )$ is a semimartingale with respect to an increasing filtration ${\mathcal P}_{t}$,
$t\in [0,T]$ and with values in a manifold $M$, we consider the process $\eta$ defined by the Stratonovich integral 
\begin{eqnarray*}
\eta (t) :=\int_0^t  t^{\xi}_{0\leftarrow s} o d\xi (s)
\end{eqnarray*}
This is a semimartingale taking values in $T_{\xi (0)} (M)$. We consider its (generalised) right-hand
time derivative (or drift) by taking conditional expectations:
\begin{eqnarray*}
D_t \eta (t)=\lim_{\epsilon \rightarrow 0}E^{{\mathcal P}_t} [\frac{\eta (t+\epsilon )-\eta (t)}{\epsilon}]
\end{eqnarray*}
Notice that if $\xi$ is a differentiable deterministic path, this notion of derivative reduces to the usual one.

Then we define the \it generalised (forward) derivative \rm
\begin{eqnarray}\label{eq12}
D^{\nabla}_t \xi (t):= t_{t\leftarrow 0} D_t \eta (t)
\end{eqnarray}
We use the symbol $\nabla$ to stress that the derivative depends on the choice of covariant derivative used to define the parallel transport, although in this work we are only consider the Levi-Civita covariant derivative.

For a (possibly time dependent)  vector field $Z$ computed along a semimartingale $\xi$, the generalised derivative is defined as

$$D_t^{\nabla} Z (t)= \lim_{\epsilon \rightarrow 0}E^{{\mathcal P}_t} [t_{t \leftarrow t+\epsilon } Z(t+\epsilon,\xi (t+\epsilon ))-Z(t,\xi (t)]$$

Let us consider our base manifold $M$ and, for a $M$-valued semimartingale $\xi$, define the
corresponding kinetic energy by
\begin{eqnarray}\label{eq13}
{\mathcal E}(\xi )=E \int_0^T ||D^{\nabla} \xi (t)||^2 dt
\end{eqnarray}

Next Theorem characterises the critical paths of $\mathcal E$. Allowed variations are processes of bounded variation $h$  satisfying $h(0)=h(T)=0$.
We have the following result:
\vskip 5mm

\begin{theorem}\label{teo1}
A diffusion process $m(\cdot )$ with generator $L_u$, $u\in C^2 (M)$, is a critical path for the energy functional ${\mathcal E}$ if and only if $D^{\nabla}_t u (t, m(t))=0$ almost everywhere or, equivalently,
\begin{eqnarray}\label{eq14}
\frac{\partial}{\partial t} u+  (\nabla_u u )+\frac{1}{2} [(\Delta u ) +\hbox{Ricci} (u)]=0
\end{eqnarray}
\end{theorem}

Notice that, in particular, we obtain the expression derived in \cite{Zambrini1999} using local coordinates.

It is shown in \cite{AZ2017} (c.f, more generally, \cite{Zambrini1999})  that the symmetries of the critical diffusion coincide with the regularisation of its classical counterpart.
In other words, if the diffusion coefficient in  (8), regarded now as variable, tends to zero, $D^\nabla \xi$ in (12)  reduces to an ordinary (strong) derivative, Eq. (14) and  the 
symmetries of the critical diffusion reduce to those of the classical functional (4).

\begin{proof}
We first write the energy functional via the lift of the process to the frame bundle, as explained in the last paragraph:
\begin{eqnarray*}
{\mathcal E}= E\int_0^T ||D_t \pi (r_x (t))||^2 dt
\end{eqnarray*}
where $D_t$ refers to the generalised derivative for processes defined on the flat space (of the Brownian motion $x$).
Then we perform variations of the Brownian motion $x$ along directions $h(\cdot )$, processes of bounded variation with $h(0)=h(T)=0$. Using  Proposition \ref{prop1}, these 
variations will give
rise to variations on the path space of the manifold $M$ along semimartingales $\zeta (\cdot )$,
where $\zeta$ is given by (\ref{eq11}).
More precisely we have,
\begin{eqnarray*}
&\frac{d}{d\epsilon}_{|_{\epsilon =0}} E\int_0^T  ||D_t \pi (r_{x+\epsilon h} (t))||^2
dt=2E\int_0^T <D_t \pi (r_x (t)), D_t \pi^{\prime} (\frac{d}{d\epsilon}_{|_{\epsilon =0}}
r_{x+\epsilon h} (t))>dt\\
&=2E\int_0^T <D_t \pi (r_x (t)), D_t  (\zeta )(t))>dt\\
&=2E\int_0^T <D_t \pi (r_x (t)), \dot h -\frac{1}{2}\hbox{Ricci} (h )-{\mathcal D}(h)(t)>dt
\end{eqnarray*}

Using integration by parts in time, the assumption
$h(0)=h(t)=0$ and the fact that there is no It\^o's extra term in the integration since $h$ is of bounded variation,
 the first term is equal to
$-2E\int_0^T  <D_t D_t \pi (r_x), h(t)>$.

We arrive to the conclusion that a process $r_x$  of the form (\ref{eq9}) is critical for the 
action functional $\mathcal E$ if and only if 
$D_t^{\nabla} u=0$
almost everywhere, which proves the Theorem.
\end{proof}

\section{Stochastic geodesics on Lie groups}
\label{sec:6}

Let $G$ denote a Lie group endowed with a left invariant metric $<~>$ and a left invariant connection $\nabla$, that we assume here to be the Levi-Civita connection. The corresponding Lie algebra $\mathcal G$ can be identified with the
tangent space $T_e G$, where $e$ is the identity element of the group. Taking a sequence of
vectors $H_k \in {\mathcal G}$, consider the following Stratonovich stochastic differential equation on the group:
\begin{eqnarray}\label{eq15}
dg(t) =T_e L_{g(t)} \left(\sum_k H_k \circ dx^k (t) -\frac{1}{2} \nabla_{H_k} H_k dt +u(t)dt \right)
\end{eqnarray}
with $g(0=e$, where $T_a L_{g(t)}:T_a G \rightarrow T_{g(t)a} G$ is the differential of the left translation $L_{g(t)}(x)=g(t)x, x\in G$ and where $x^k (t)$ are independent real valued Brownian motions. The vector $u(\cdot )$ is assumed to be non random, $u(\cdot )\in C^2 ([0,T]; {\mathcal G})$.

The stochastic energy functional for a general $G$-vaued semimartingale $\xi (t), t\in [0,T]$, reads:
\begin{eqnarray}\label{eq16}
{\mathcal E} (\xi )= E\int_0^T ||T_{\xi (t)} L_{\xi (t)^{-1}} D^{\nabla}_t \xi (t)||^2 dt
\end{eqnarray}
Assume furthermore that $\nabla_{H_k} H_k =0$ for all $k$ (in particular the Stratanovich integral in (\ref{eq15}) coincides with the It\^o one). Then the following result holds:

\begin{theorem}[\cite{ACC2014}]
A $G$-valued semimartingale of the form (\ref{eq15}) is critical for the energy functional (\ref{eq16}) if and only if 
the vector field $u(\cdot )$ satisfies the equation
\begin{eqnarray*}
\frac{d}{dt} u(t)= \hbox{ad}_{u (t)} u(t) -\frac{1}{2} \left( \sum_k \nabla_{H_k} \nabla_{H_k} u (t)
+\hbox{Ricci} (u(t))\right)
\end{eqnarray*}
\end{theorem}

When $H_k=0$ for all $k$ the equation reduces to the well known Euler-Poincar\'e equation
for (deterministic) geodesics in Lie groups $\frac{d}{dt} u(t)= \hbox{ad}_{u (t)} u(t)$.

Up to some sign changes, the right invariant case is analogous.

The theorem also holds for infinite-dimensional Lie groups and allows, as a particular case, to derive the Navier-Stokes equation, when
the problem is formulated on the diffeomorphisms group (c.f. \cite{ACC2014}).

\section{Relation with stochastic forward-backward differential equations}
\label{sec:7}

Deterministic geodesics solve second order differential equations and as such can be obtained
using standard methods for such equations, with given initial position and velocity as well as
with initial and final given positions. The meaning of "second order" stochastic differential equations is not so clear. A possible method is its characterisation via stochastic 
forward-backward differential equations. In local coordinates (c.f. notations defined in (\ref{eq8}), a stochastic geodesic in the time interval $[0,T]$ reads

\begin{eqnarray}
m^j (t)= m^j (0) +\int_0^t \sigma_k^j (m(s))dx^k (s)-\int_0^t \left( \frac{1}{2} g^{m,n}\Gamma_{m,n}^j  (m(s))-y^j (s) \right) ds\nonumber \\
y^j (t)=y^j (T) -\int_t^T Z_k^j (s) dx^k (s)-\frac{1}{2}\int_t^T \hbox{Ricci}^j (m(s))ds\nonumber
\end{eqnarray}

Given $m^j (0)$ and $y(T)=u(T, m(T))$ these kind of systems may provide solutions of the form
$(m(t), y(t))$ with $y(t)=u(t, m(t))$ corresponding to our stochastic geodesics (c.f., for example, 
\cite{Delarue2002}). The term $Z$ is an a priori unknown of the equation, but is in fact determined a posteriori by the solution $(m,y)$.

In the case of stochastic geodesics on Lie groups, the characterisation via forward-backward equations was described in \cite{CC2013}.
An extension to  infinite dimensional Lie groups and, in particular, to the Navier-Stokes equation framework, is also possible (\cite{CS2009}, \cite{CCQ2013}).

\vskip 10mm

\noindent{\bf Acknowledgements}\\

The authors acknowledge the support of the FCT Portuguese grant PTDC/MAT-STA/28812/2017.

\end{document}